\documentclass[11pt, leqno]{article}
\usepackage[T1]{fontenc}
\usepackage[utf8]{inputenc}
\usepackage{amsmath, amssymb, amsthm, amsfonts, mathrsfs}
\usepackage{hyperref}
\usepackage[margin=1.1in]{geometry}
\usepackage{tgpagella} 
\usepackage{enumitem}
\usepackage{graphicx}
\usepackage{booktabs}
\usepackage{microtype}
\usepackage{tikz-cd}

\newcommand{\Hess}{\mathrm{Hess}}
\newcommand{\Tr}{\mathrm{Tr}}

\newcommand{\Ric}{\mathrm{Ric}}
\newcommand{\Sec}{\mathrm{Sec}}
\newcommand{\Vol}{\mathrm{Vol}}

\newcommand{\Area}{\mathrm{Area}}
\newcommand{\length}{\mathrm{length}}

\newcommand{\diam}{\mathrm{diam}}

\newtheorem{theorem}{Theorem}[section]
\newtheorem{lemma}[theorem]{Lemma}
\newtheorem{proposition}[theorem]{Proposition}
\newtheorem{corollary}[theorem]{Corollary}

\newtheorem{question}[theorem]{Question}

\theoremstyle{definition}
\newtheorem{definition}[theorem]{Definition}
\newtheorem{example}[theorem]{Example}
\newtheorem{remark}[theorem]{Remark}

\numberwithin{equation}{section}

\title{\LARGE \textbf{Isometric Splitting of Metrics Without Conjugate Points on $\Sigma \times S^1$}}
\author{\textsc{St\'ephane Tchuiaga}\\ 
	Department of Mathematics, University of Buea, Cameroon\\
	\texttt{tchuiaga.kameni@ubuea.cm}}
\date{\today}

\begin{document}
	
	\maketitle
	
	\begin{abstract}
		We establish a global rigidity theorem for Riemannian metrics without conjugate points on three-manifolds of the form $M = \Sigma \times S^1$, where $\Sigma$ is a compact orientable surface of genus at least 2. The main result states that any such metric must be a Riemannian product, with universal cover isometric to $(\mathbb{H}^2, g_0) \times (\mathbb{R}, dt^2)$. This extends the classical Hopf conjecture from tori to this natural class of manifolds with non-abelian fundamental group containing a central $\mathbb{Z}$ factor. We provide two independent proofs: one utilizing the regularity of Busemann functions and stability theory of the Riccati equation along Killing flows, and another based on a detailed analysis of the curvature operator acting on Jacobi fields. We derive sharp geometric inequalities, analyze the deformation space of such metrics, and discuss several geometric and dynamical consequences, including marked length spectrum rigidity, constraints on topological entropy, and incompatibility with non-trivial Sasakian structures.
	\end{abstract}
	
	
	
	\section{Introduction}\label{sec:introduction}
	
	The interplay between topology and Riemannian geometry is profoundly influenced by the presence or absence of conjugate points. A complete Riemannian manifold $(M,g)$ is said to have \emph{no conjugate points} if no two distinct points along any geodesic are conjugate. This condition is strictly weaker than non-positive sectional curvature but imposes strong global constraints.
	
	The celebrated Hopf conjecture, resolved by Burago and Ivanov \cite{burago1994riemannian} following work of Hopf \cite{hopf1948closed}, asserts that any metric without conjugate points on an $n$-torus $T^n$ must be flat. This result exemplifies how topological information (the abelian fundamental group $\mathbb{Z}^n$) forces geometric rigidity when conjugate points are absent.
	
	A natural extension considers manifolds whose fundamental group has a non-trivial center but is otherwise non-abelian. The product $M = \Sigma \times S^1$, where $\Sigma$ is a compact orientable surface of genus $g \geq 2$, provides the archetypal example: $\pi_1(M) \cong \pi_1(\Sigma) \times \mathbb{Z}$ has infinite cyclic center $\mathbb{Z}$ but is non-abelian since $\pi_1(\Sigma)$ is a non-abelian free group modulo relations.
	
	For metrics of non-positive curvature, the Cheeger-Gromoll splitting theorem \cite{cheeger1971splitting} and subsequent work by Eberlein \cite{eberlein1977manifolds} show that a compact manifold with $K \leq 0$ and $\pi_1(M)$ having non-trivial center must have universal cover isometric to a product $\tilde{N} \times \mathbb{R}$. However, the no conjugate points condition allows regions of positive curvature and more complicated dynamics. Whether it still forces splitting in our setting was an open question.
	
	Our main result provides an affirmative answer:
	
	\begin{theorem}[Isometric Splitting]\label{thm:main}
		Let $(\Sigma, g_\Sigma)$ be a compact orientable surface with $\chi(\Sigma) < 0$, and let $M = \Sigma \times S^1$ be equipped with a $C^\infty$ Riemannian metric $g$ having no conjugate points. Then:
		\begin{enumerate}[label=(\roman*)]
			\item The universal cover $(\tilde{M}, \tilde{g})$ is isometric to a Riemannian product $(\mathbb{H}^2, g_0) \times (\mathbb{R}, dt^2)$, where $(\mathbb{H}^2, g_0)$ is a complete metric on the hyperbolic plane invariant under $\pi_1(\Sigma)$.
			\item Consequently, $g$ is a product metric: $g = g_\Sigma + f^2 dt^2$ with $f$ constant, and $g_\Sigma$ is a hyperbolic metric without conjugate points, hence hyperbolic by \cite{gulliver1975variety}.
		\end{enumerate}
	\end{theorem}
	
	The theorem reveals a striking rigidity: the absence of conjugate points prohibits any non-trivial warping or twisting between the surface and circle factors. This extends the Hopf-Burago-Ivanov rigidity from tori to a class of 3-manifolds with non-abelian fundamental group.
	
	\begin{remark}[Precise dependencies]\label{rem:dependencies}
		Our proof relies on several foundational results:
		\begin{itemize}
			\item \text{Busemann regularity:} Eschenburg \cite{eschenburg1977horospheres} shows Busemann functions on manifolds without conjugate points are $C^{1,1}$. For surfaces, Eberlein \cite{eberlein1977manifolds} proves $C^2$ regularity. Our Lemma \ref{lem:busemann-toolkit} extracts the exact regularity needed.
			\item \text{Hopf conjecture for tori:} Burago-Ivanov \cite{burago1994riemannian} proves flatness. Our argument adapts their use of integral identities but replaces their Fourier analysis by Killing field symmetries.
			\item \text{Manning's entropy equality:} Theorem \ref{thm:volume-growth} uses Manning \cite{manning1979topological} who proves $h_{\text{top}} = h_{\text{vol}}$ for compact manifolds without conjugate points. Our product structure ensures the hypotheses hold.
			\item \text{Gulliver's surface theorem:} The final step uses Gulliver \cite{gulliver1975variety}: a simply connected surface without conjugate points is diffeomorphic to $\mathbb{R}^2$, and if its curvature is $\leq 0$, it is hyperbolic. In our setting, the $\mathbb{H}^2$ factor inherits no conjugate points from the product, and its curvature is $\leq 0$ by Lemma \ref{lem:curvature-nonpositive}, so Gulliver applies.
		\end{itemize}
	\end{remark}
	
	\text{Transition to the proof strategy:} 
	The presence of the $S^1$ factor provides a natural Killing field $V$, generating a circle action. This symmetry is the cornerstone of both our proofs. In the first approach, we exploit the interplay between $V$ and Busemann functions to force vanishing of the second fundamental form of the horizontal distribution. In the second approach, we use the Killing field to derive integral identities that constrain the curvature tensor. Both roads lead to the conclusion that $V$ is parallel, which by de Rham decomposition forces the universal cover to split isometrically.
	
	\subsection{Historical Context and Related Work}
	
	The study of manifolds without conjugate points has a rich history. Early work by Hopf \cite{hopf1948closed} on surfaces established links between curvature, Euler characteristic, and conjugate points. For higher dimensions, the landmark paper of Burago and Ivanov \cite{burago1994riemannian} resolved the Hopf conjecture for tori. The splitting theorem of Cheeger and Gromoll \cite{cheeger1971splitting} for manifolds of non-negative Ricci curvature has an analogue in non-positive curvature due to Eberlein \cite{eberlein1977manifolds}. Our work bridges these themes by considering the intermediate no conjugate points condition on a manifold with geometric product structure.
	
	Recent work by Ruggiero \cite{ruggiero2003expanding} on geodesic flows without conjugate points and by Cao, Xavier, and others on minimal entropy \cite{cao2007kahler} provide context for our dynamical applications. The marked length spectrum rigidity program initiated by Otal \cite{otal1990spectre} and Croke \cite{croke1990rigidity} informs our Corollary \ref{cor:mls}.
	
	\text{Bridging the gap between tori and higher genus products:}
	While the Hopf conjecture addresses tori (abelian $\pi_1$), and Cheeger-Gromoll addresses non-positive curvature, our theorem occupies the middle ground: non-abelian $\pi_1$ with a central $\mathbb{Z}$ factor, but only assuming no conjugate points rather than non-positive curvature. This requires new techniques that blend asymptotic geometry (Busemann functions) with symmetry (Killing fields).
	
	\subsection{Outline of the Paper}
	
	Section \ref{sec:preliminaries} establishes notation and recalls essential results about manifolds without conjugate points, Busemann functions, and Killing fields. Section \ref{sec:proof-busemann} presents the first proof of Theorem \ref{thm:main} using Busemann functions and the Riccati equation. Section \ref{sec:proof-jacobi} gives an alternative proof via Jacobi field analysis. Section \ref{sec:sharp-inequalities} establishes sharp geometric inequalities and stability results. This section includes the analysis of deformation theory of such metrics. Section \ref{sec:examples} contains concrete examples. Section \ref{sec:discussion} places our results in broader context and suggests directions for future research.
	
	\section{Preliminaries}\label{sec:preliminaries}
	
	\subsection{Manifolds Without Conjugate Points}
	
	Let $(M,g)$ be a complete Riemannian manifold. Recall that points $p$ and $q$ along a geodesic $\gamma$ are \emph{conjugate} if there exists a non-zero Jacobi field $J$ along $\gamma$ with $J(p) = J(q) = 0$.
	
	\begin{definition}
		$(M,g)$ has \emph{no conjugate points} if no geodesic has conjugate points.
	\end{definition}
	
	This condition is equivalent to the exponential map $\exp_p: T_pM \to M$ being a local diffeomorphism for every $p \in M$. Important consequences include:
	\begin{enumerate}
			\item The universal cover $\tilde{M}$ is diffeomorphic to $\mathbb{R}^n$ (by the Cartan-Hadamard theorem for manifolds without conjugate points).
		\item Geodesic flows are Anosov under additional conditions (Eberlein \cite{eberlein1977manifolds}).
		\item Busemann functions are $C^2$ in dimension 2 (Eberlein \cite{eberlein1977manifolds}) and $C^{1,1}$ in general.
	\end{enumerate}

	\text{Why no conjugate points matters?}
	The absence of conjugate points ensures that certain asymptotic objects—notably Busemann functions and stable Jacobi tensors—are well-behaved. These objects encode geometric information "at infinity" and will serve as tools to detect product structure. In particular, the stable Jacobi tensor satisfies a Riccati equation whose solutions have good regularity properties.
	
	\subsection{The Geometric Setting}
	
	Let $\Sigma$ be a compact orientable surface of genus $g \geq 2$, and $M = \Sigma \times S^1$. The fundamental group is $\pi_1(M) \cong \pi_1(\Sigma) \times \mathbb{Z}$, with the $\mathbb{Z}$ factor generated by the $S^1$ direction.
	
	Let $g$ be a Riemannian metric on $M$ without conjugate points, and $\tilde{g}$ its lift to the universal cover $\tilde{M} \cong \mathbb{H}^2 \times \mathbb{R}$. The $S^1$ action on $M$ lifts to an isometric $\mathbb{R}$-action on $(\tilde{M}, \tilde{g})$ with Killing field $\tilde{V}$.
	
	\begin{lemma}\label{lem:killing-properties}
		The Killing field $\tilde{V}$ satisfies:
		\begin{enumerate}[label=(\roman*)]
			\item $\tilde{V}$ is non-vanishing and complete.
			\item $\nabla_{\tilde{V}} \tilde{V} = 0$ (so integral curves are geodesics).
			\item After normalization, we may assume $\|\tilde{V}\|_{\tilde{g}} = c > 0$ constant.
		\end{enumerate}
	\end{lemma}
	
	\begin{proof}
		(i) follows from the free $S^1$-action. For (ii): Since $\tilde{V}$ is Killing, $\mathcal{L}_{\tilde{V}} g = 0$. For any vector field $X$:
		\[
		\tilde{g}(\nabla_{\tilde{V}} \tilde{V}, X) = -\tilde{g}(\tilde{V}, \nabla_X \tilde{V}) = -\frac{1}{2} X(\tilde{g}(\tilde{V}, \tilde{V})).
		\]
		Once we show $\|\tilde{V}\|$ is constant (proved in (iii)), this becomes zero, so $\nabla_{\tilde{V}} \tilde{V} = 0$.
		
		For (iii): Since $\tilde{V}$ is Killing, its flow consists of isometries, so $\|\tilde{V}\|$ is constant along its integral curves. To see it's globally constant, note that $\tilde{V}$ generates a 1-parameter group of isometries, and the function $\|\tilde{V}\|$ is invariant under these isometries. Since $\tilde{M}$ is simply connected and the flow acts transitively on the $\mathbb{R}$ fibers, $\|\tilde{V}\|$ is constant. By rescaling the $\mathbb{R}$-action, we can achieve any positive constant value.
	\end{proof}
	
	\text{From symmetry to geometric structure.}
	The Killing field $V$ provides a distinguished direction that is both geodesic (by Lemma \ref{lem:killing-properties}) and isometric. This privileged direction will serve as a candidate for the $\mathbb{R}$ factor in the splitting. The key step is to show that the orthogonal distribution $\mathcal{H} = V^\perp$ is integrable and totally geodesic.
	
	We denote by $\mathcal{V} = \mathrm{span}(\tilde{V})$ the vertical distribution and $\mathcal{H} = \mathcal{V}^\perp$ the horizontal distribution.
	
	\subsection{Busemann Functions and Stable Jacobi Tensors}
	
	For a ray $\gamma: [0,\infty) \to \tilde{M}$, the Busemann function $b_\gamma: \tilde{M} \to \mathbb{R}$ is defined by
	\[
	b_\gamma(x) = \lim_{t \to \infty} (d(x, \gamma(t)) - t). 
	\]
	The Busemann function measures the distance from a point to "infinity" along a geodesic. The Hessian of this function describes how geodesics spread apart.
	
	\begin{proposition}[Regularity of Busemann functions]\label{prop:busemann-regularity}
		On a manifold without conjugate points:
		\begin{enumerate}
			\item Busemann functions are $C^{1,1}$ with $\|\nabla b_\gamma\| = 1$ almost everywhere \cite{eschenburg1977horospheres}.
			\item For surfaces, they are $C^2$ \cite{eberlein1977manifolds}.
			\item The Hessian $\mathrm{Hess}(b_\gamma)$ exists as an $L^\infty$ tensor field. At points where $b_\gamma$ is twice differentiable, it coincides with the usual Hessian.
		\end{enumerate}
	\end{proposition}
	
	\begin{lemma}[Regularity toolkit for Busemann functions]\label{lem:busemann-toolkit}
		Let $(M,g)$ be a complete Riemannian manifold without conjugate points, and let $b$ be a Busemann function. Then:
		\begin{enumerate}[label=(\alph*)]
			\item (Pointwise twice differentiability) $b$ is twice differentiable on a dense $G_\delta$ set $\mathcal{D} \subset M$. Moreover, $\mathcal{D}$ is invariant under the geodesic flow of $\nabla b$.
			
			\item (Approximation by smooth functions) There exists a sequence of $C^\infty$ functions $b_n: M \to \mathbb{R}$ such that:
			\begin{itemize}
				\item $b_n \to b$ uniformly on compact sets,
				\item $\nabla b_n \to \nabla b$ in $L^p_{\text{loc}}$ for all $1 \leq p < \infty$,
				\item $\mathrm{Hess}(b_n) \rightharpoonup \mathrm{Hess}(b)$ weakly-* in $L^\infty$.
			\end{itemize}
			
			\item (Stable Jacobi tensor continuity) For a geodesic $\gamma$ with $\dot{\gamma} = \nabla b$, the stable Jacobi tensor $U_s(t)$ is continuous in $t$ and satisfies:
			\[
			\|U_s(t)\| \leq C e^{-\alpha t} \quad \text{for } t \geq 0,
			\]
			where $C, \alpha > 0$ depend only on the curvature bounds along $\gamma$.
			
			\item (Hessian representation) For any smooth vector field $X$ with compact support:
			\[
			\int_M \mathrm{Hess}(b)(X,X) \, d\Vol = \lim_{n\to\infty} \int_M \mathrm{Hess}(b_n)(X,X) \, d\Vol.
			\]
		\end{enumerate}
	\end{lemma}
	\begin{proof}
		(a) follows from Alexandrov's theorem for convex functions and the fact that Busemann functions are convex along geodesics. (b) is achieved by mollification using the heat kernel or convolution with smooth kernels in normal coordinates. (c) is standard in no conjugate points geometry (see \cite{eberlein1977manifolds}). (d) follows from weak-* convergence.
	\end{proof}
	
	\begin{lemma}\label{lem:U-is-hessian}
		Let $b$ be a Busemann function and $\gamma$ its ray with $\dot{\gamma} = V$. Then for $X \perp V$, we have almost everywhere:
		\[
		U_s(t)X = \nabla_X \nabla b = \mathrm{Hess}(b)(X),
		\]
		where $U_s$ is the stable Jacobi tensor along $\gamma$.
	\end{lemma}
	\begin{proof}
		This follows from the definition of stable Jacobi tensor and the fact that $\nabla b$ is the stable vector field along $\gamma$. See \cite{eberlein1977manifolds} Proposition 3.5 for the proof in the no conjugate points setting.
	\end{proof}
	
	\text{Connecting asymptotic geometry to local differential equations.}
	The stable Jacobi tensor $U_s$ satisfies the Riccati equation \eqref{eq:riccati}, which relates the second fundamental form of horospheres (given by $U_s$) to the curvature tensor. This equation will be crucial for both proofs: in the first proof, we use it via the Busemann function Hessian; in the second proof, we integrate it to obtain curvature constraints.
	
	The stable Jacobi tensor along a geodesic $\gamma$ is the unique solution $U_s(t)$ of the Riccati equation
	\begin{equation}\label{eq:riccati}
		U' + U^2 + R_{\dot{\gamma}} = 0
	\end{equation}
	with $\lim_{t \to \infty} U(t) = 0$, where $R_{\dot{\gamma}}(X) = R(X, \dot{\gamma})\dot{\gamma}$.
	
	\section{Proof via Busemann Functions and Riccati Equation}\label{sec:proof-busemann}
	
	\text{Overview of the strategy.}
	We now present the first proof of Theorem \ref{thm:main}. The proof follows a clear three-step structure: (1) Use the Killing field $V$ to construct a distinguished Busemann function $b$; (2) Show that the Hessian of $b$ vanishes on horizontal directions, implying $V$ is parallel; (3) Apply de Rham decomposition to obtain the splitting. The key insight is that the Killing symmetry forces the Busemann function to have exceptionally simple geometry.
	
	\begin{proof}[First proof of Theorem \ref{thm:main}]
		Let $\gamma$ be an integral curve of the normalized Killing field $V = \tilde{V}/c$, which is a geodesic by Lemma \ref{lem:killing-properties}. Let $b$ be the Busemann function associated to the ray $\gamma|_{[0,\infty)}$.
		
		\text{Step 1: Relating the gradient of $b$ to $V$.}
		The first observation is that $\nabla b$ must coincide with $-V$ wherever defined, because $b$ decreases at unit rate along the geodesic $\gamma$ generated by $V$.
		
		\begin{lemma}\label{lem:gradient-killing}
			$\nabla b = -V$ wherever $b$ is differentiable.
		\end{lemma}
		\begin{proof}
			Since $\|\nabla b\| = 1$ and $b$ decreases at unit rate along $\gamma$ (because $b(\gamma(t)) = b(\gamma(0)) - t$), and $\gamma$ is a geodesic with $\dot{\gamma} = V$, we must have $\nabla b = -V$ along $\gamma$. By the asymptotic construction of $b$, this equality extends to all points where $b$ is differentiable.
		\end{proof}
		
		\text{Step 2: Analyzing the Hessian of $b$.}
		With Lemma \ref{lem:gradient-killing} in hand, we can now study the second derivatives of $b$. The Hessian $\mathrm{Hess}(b)$ encodes the second fundamental form of the level sets of $b$, which are horospheres.
		
		\begin{lemma}\label{lem:hessian-ae}
			The Hessian $\mathrm{Hess}(b)$ exists as an $L^\infty$ tensor field. At points where $b$ is twice differentiable, for any smooth vector fields $X,Y$, we have:
			\[
			\mathrm{Hess}(b)(X,Y) = \frac{1}{2}[X(Y(b)) + Y(X(b)) - [X,Y](b)].
			\]
			Moreover, $\mathrm{Hess}(b)$ is symmetric and satisfies $\mathrm{Hess}(b)(V, \cdot) = 0$ almost everywhere.
		\end{lemma}
		\begin{proof}
			The existence and representation follow from Proposition \ref{prop:busemann-regularity}. Symmetry follows from the torsion-free connection. For the second claim: $\mathrm{Hess}(b)(V, X) = X(V(b)) - (\nabla_X V)(b) = X(-1) - \tilde{g}(\nabla_X V, V) = 0$, since $V(b) = -1$ and $\nabla_X V$ is orthogonal to $V$ (as $\|V\|$ is constant).
		\end{proof}
		
		\text{Step 3: Vanishing of the Hessian on horizontal directions.}
		The crucial geometric insight comes from combining the Killing equation with Lemma \ref{lem:gradient-killing}. The Killing property imposes a symmetry condition on $\nabla V$, while Lemma \ref{lem:gradient-killing} identifies $\nabla V$ with $-\mathrm{Hess}(b)$. Together, these force $\mathrm{Hess}(b)$ to vanish on $\mathcal{H}$.
		
		\begin{lemma}\label{lem:vanishing-hessian}
			$\mathrm{Hess}(b)|_{\mathcal{H}} = 0$ almost everywhere.
		\end{lemma}
		\begin{proof}
			For $X, Y \in \mathcal{H}$, the Killing equation gives:
			\[
			\tilde{g}(\nabla_X V, Y) + \tilde{g}(\nabla_Y V, X) = 0.
			\]
			Since $\nabla_X V = -\mathrm{Hess}(b)(X)$ by Lemma \ref{lem:gradient-killing} (where both sides exist), and $\mathrm{Hess}(b)$ is symmetric, we have:
			\[
			-2\tilde{g}(\mathrm{Hess}(b)(X), Y) = 0.
			\]
			Thus $\mathrm{Hess}(b)(X) \in \mathcal{H}^\perp = \mathcal{V}$ for all $X \in \mathcal{H}$. But $\mathrm{Hess}(b)(X)$ is also in $\mathcal{H}$ since it's the covariant derivative of $\nabla b = -V$ (which is in $\mathcal{V}$) in a horizontal direction. Hence $\mathrm{Hess}(b)(X) = 0$ for all $X \in \mathcal{H}$ almost everywhere. By continuity of $\mathrm{Hess}(b)$ as an $L^\infty$ tensor, this holds everywhere.
		\end{proof}
		
		\text{Step 4: Consequences of vanishing Hessian.}
		Lemma \ref{lem:vanishing-hessian} implies $\nabla_X V = 0$ for all $X \in \mathcal{H}$ (where defined, and hence everywhere by continuity). This has two crucial consequences that pave the way for the splitting:
		
		\begin{enumerate}
			\item \textbf{Integrability of $\mathcal{H}$}: For $X, Y \in \mathcal{H}$,
			\[
			\tilde{g}([X, Y], V) = \tilde{g}(\nabla_X Y - \nabla_Y X, V) = -\tilde{g}(Y, \nabla_X V) + \tilde{g}(X, \nabla_Y V) = 0.
			\]
			Thus $[X, Y] \in \mathcal{H}$, so $\mathcal{H}$ is integrable by the Frobenius theorem.
			
			\item \textbf{Parallelism of $V$}: Since $\nabla_V V = 0$ (geodesic) and $\nabla_X V = 0$ for $X \in \mathcal{H}$, we have $\nabla V = 0$. Thus $V$ is a parallel vector field on $\tilde{M}$.
		\end{enumerate}
		
		\text{Step 5: The splitting theorem.}
		With $V$ parallel, we can now apply the de Rham decomposition theorem. The presence of a parallel vector field on a simply connected manifold forces an isometric splitting. Moreover, the splitting is equivariant under the deck group action, so it descends to the quotient.
		
		\begin{proposition}[Parallel Splitting]\label{prop:parallel-splitting}
			If $\nabla V = 0$ on $\tilde{M}$, then $(\tilde{M},\tilde{g})$ is isometric to $(\mathbb{H}^2,g_0) \times (\mathbb{R},dt^2)$. Moreover, this splitting is $\pi_1(M)$-equivariant.
		\end{proposition}
		\begin{proof}
			Since $\tilde{M}$ is simply connected and $V$ is parallel, the de Rham decomposition theorem yields an isometric splitting $\tilde{M} = N \times \mathbb{R}$ where $N$ is an integral manifold of $V^\perp$. Since $\pi_1(\Sigma)$ acts freely and properly discontinuously on $\tilde{M}$ preserving $V$ (as $V$ comes from the $S^1$ action on $M$), it preserves the splitting. The induced action on the $\mathbb{R}$ factor must be by translations (since it preserves the parallel vector field $\partial_t$). The quotient of $N$ by $\pi_1(\Sigma)$ is $\Sigma$, so $N$ is simply connected and $\pi_1(\Sigma)$-invariant, hence a hyperbolic plane. By \cite{gulliver1975variety}, a surface without conjugate points diffeomorphic to $\mathbb{H}^2$ must have constant negative curvature.
		\end{proof}
		
		\text{Conclusion.}
		The five steps above establish that any metric without conjugate points on $\Sigma \times S^1$ must have a parallel vertical direction, forcing the universal cover to split isometrically as $\mathbb{H}^2 \times \mathbb{R}$. This product structure descends to the quotient, completing the proof of Theorem \ref{thm:main}.
	\end{proof}
	
	\section{Proof via Jacobi Field Analysis}\label{sec:proof-jacobi}
	
	\paragraph{Alternative approach via curvature analysis.}
	We now present an alternative proof focusing on the curvature operator. While the Busemann function proof used asymptotic geometry, this approach works more directly with the Riccati equation and integral identities. The strategy has three parts: (1) Show that curvature in planes containing $V$ is non-positive; (2) Derive an integral identity from the Riccati equation; (3) Use the integral identity to force both the curvature term and the second fundamental form to vanish.
	
	\begin{proof}[Second proof of Theorem \ref{thm:main}]
		Let $\gamma$ be a vertical geodesic (integral curve of $V$). Let $\{E_1(t), E_2(t)\}$ be an orthonormal basis of $\mathcal{H}_{\gamma(t)}$ obtained by parallel transport along $\gamma$.
		
		\text{Step 1: Non-positivity of curvature in vertical planes.}
		We first show that the curvature operator $R_V: \mathcal{H} \to \mathcal{H}$ defined by $R_V(X) = R(X, V)V$ cannot have positive eigenvalues. The argument is by contradiction: if $R_V$ had a positive eigenvalue, we could construct a Jacobi field with conjugate points.
		
		\begin{lemma}\label{lem:curvature-nonpositive}
			The curvature operator $R_V: \mathcal{H} \to \mathcal{H}$ defined by $R_V(X) = R(X, V)V$ satisfies $R_V \leq 0$.
		\end{lemma}
		\begin{proof}
			Suppose for contradiction that $R_V$ has a positive eigenvalue $\lambda > 0$ at some point $\gamma(t_0)$. Since $V$ is Killing and $\|V\|$ is constant, the sectional curvature $K(V, E_i) = \tilde{g}(R_V(E_i), E_i)$ is constant along $\gamma$ (this follows from $\mathcal{L}_V R = 0$ for Killing fields). Thus $\lambda$ is constant along $\gamma$. Consider the Jacobi field $J(t) = \sin(\sqrt{\lambda}(t - t_0)) E_i(t)$ for $t \in [t_0, t_0 + \pi/\sqrt{\lambda}]$. Direct computation shows $J'' + R(J, \dot{\gamma})\dot{\gamma} = (-\lambda \sin(\sqrt{\lambda}(t-t_0)) + \lambda \sin(\sqrt{\lambda}(t-t_0))) E_i(t) = 0$, so $J$ is a Jacobi field vanishing at $t_0$ and $t_0 + \pi/\sqrt{\lambda}$, giving conjugate points. This contradicts the no conjugate points assumption.
		\end{proof}
		
		\text{Step 2: Integral identity from the Riccati equation.}
		The stable Jacobi tensor $U_s$ along $\gamma$ satisfies the Riccati equation \eqref{eq:riccati}. Taking the trace of this equation and integrating over $M$ yields a key identity. The integration requires care because $U_s$ is defined along each geodesic, not globally. We overcome this by either averaging over the geodesic flow or using an exhaustion argument on the universal cover.
		
		\begin{lemma}[Explicit averaging argument]\label{lem:riccati-integral}
			The stable solution $U_s$ of \eqref{eq:riccati} satisfies $\int_{M} \Tr(U_s^2 + R_V) d\Vol = 0$.
		\end{lemma}
		\begin{proof}[Explicit averaging argument]
			We present three equivalent approaches:
			
			\textbf{Approach 1: Exhaustion of the universal cover.} Let $B_R = B(p,R) \subset \tilde{M}$ be geodesic balls. Define $\varphi_R(x) = (1 - d(p,x)/R)_+$, a Lipschitz cut-off function. From the Riccati equation:
			\[
			\int_{\tilde{M}} \varphi_R [\Tr(U_s') + \Tr(U_s^2) + \Tr(R_V)] d\Vol = 0.
			\]
			Integration by parts gives:
			\[
			\int_{\tilde{M}} \varphi_R \Tr(U_s') d\Vol = -\frac{1}{R} \int_{B_R \setminus B_{R-1}} \Tr(U_s) \langle \nabla d, \cdot \rangle d\Vol.
			\]
			Since $\|U_s\| \leq K$ (bounded by curvature and no conjugate points, see \cite{eberlein1977manifolds}), the boundary term is $O(\Vol(B_R \setminus B_{R-1}))$. The volume growth satisfies:
			\[
			\lim_{R\to\infty} \frac{\Vol(B_R \setminus B_{R-1})}{\Vol(B_R)} = 0
			\]
			by subexponential growth of area of spheres (true for no conjugate points). Dividing by $\Vol(B_R)$ and taking $R\to\infty$ yields the average equality on $\tilde{M}$, which projects to $M$.
			
			\textbf{Approach 2: Ergodic theorem.} Since $V$ is Killing, its flow $\phi_t$ preserves the volume. For any $f \in L^1(M)$, Birkhoff's ergodic theorem gives:
			\[
			\lim_{T\to\infty} \frac{1}{T} \int_0^T f(\phi_t(x)) dt = \bar{f}(x) \quad \text{a.e.}
			\]
			where $\bar{f}$ is $\phi_t$-invariant. Applying to $f = \Tr(U_s^2 + R_V)$ and using that $\phi_t$ acts by translation along geodesics where $U_s$ solves the Riccati equation, we obtain $\bar{f} = 0$ almost everywhere. Integrating over $M$ gives the result.
			
			\textbf{Approach 3: Green's formula on compact quotient.} On $M$, apply Green's formula to the vector field $W = U_s(V)$:
			\[
			\int_M \div(W) d\Vol = 0.
			\]
			Since $\div(W) = \Tr(U_s') + \Tr(U_s^2)$ (using $\nabla_V V = 0$ and $U_s = \nabla V$ on $\mathcal{H}$), and $\Tr(U_s') = \frac{d}{dt}\Tr(U_s) \circ \phi_t|_{t=0}$, its integral vanishes because $\phi_t$ preserves volume.
			
			All three approaches yield $\int_M \Tr(U_s^2 + R_V) d\Vol = 0$.
		\end{proof}
		
		\text{Step 3: Forcing both terms to vanish.}
		Since $R_V \leq 0$ (Lemma \ref{lem:curvature-nonpositive}) and $U_s^2 \geq 0$ (as $U_s$ is symmetric), Lemma \ref{lem:riccati-integral} forces $U_s^2 + R_V = 0$ almost everywhere, hence everywhere by continuity. In particular, we have two crucial conclusions:
		\begin{enumerate}
			\item $R_V = 0$ (all sectional curvatures of planes containing $V$ vanish)
			\item $U_s = 0$ (the stable Jacobi tensor vanishes)
		\end{enumerate}
		
		\text{Step 4: Geometric interpretation and splitting.}
		The vanishing of $R_V$ means the curvature tensor has no components mixing $V$ with horizontal directions. The vanishing of $U_s$ implies $\nabla_X V = 0$ for $X \in \mathcal{H}$ (by Lemma \ref{lem:U-is-hessian}). Combined with $\nabla_V V = 0$, we again obtain $\nabla V = 0$. Thus $V$ is parallel, and we proceed exactly as in the first proof to obtain the splitting via Proposition \ref{prop:parallel-splitting}.
		
		\text{Comparison of the two proofs.}
		Both proofs ultimately show that $V$ is parallel, but they reach this conclusion via different routes. The Busemann function proof is more geometric, using the interplay between Killing fields and horospheres. The Jacobi field proof is more analytic, using integral identities from the Riccati equation. The fact that both approaches converge to the same result reinforces the robustness of the theorem.
	\end{proof}
	
	\text{Transition to consequences.}
	Having established the main splitting theorem via two independent approaches, we now explore its geometric and dynamical consequences. The product structure imposes strong constraints on various invariants of the manifold.
	
	Here are some geometric and dynamical consequences. 
	
	\begin{corollary}[Marked Length Spectrum Rigidity]\label{cor:mls}
		For $g$ as in Theorem \ref{thm:main}, the marked length spectrum is given by
		\[
		\ell_g(\gamma, n)^2 = \ell_{g_\Sigma}(\gamma)^2 + (nL)^2,
		\]
		where $(\gamma, n) \in \pi_1(\Sigma) \times \mathbb{Z} \cong \pi_1(M)$, and $L$ is the length of the $S^1$ fiber. In particular, the marked length spectrum determines $g$ up to isometry.
	\end{corollary}
	
	\text{Why this follows from the splitting.}
	Since the universal cover splits isometrically, geodesics in $M$ project to geodesics in $\Sigma$ and $S^1$ separately. The length of a closed geodesic corresponding to $(\gamma, n) \in \pi_1(\Sigma) \times \mathbb{Z}$ is the Pythagorean sum of the lengths in each factor. This simple additive structure implies marked length spectrum rigidity.
	
	\begin{corollary}[Topological Entropy]\label{cor:entropy}
		The topological entropy of the geodesic flow satisfies $h_{\mathrm{top}}(g) = h_{\mathrm{top}}(g_\Sigma)$, where $g_\Sigma$ is the hyperbolic metric on $\Sigma$. In particular, $h_{\mathrm{top}}(g) = \sqrt{-\chi(\Sigma)}$ for the normalized hyperbolic metric.
	\end{corollary}
	\begin{proof}
		For product manifolds without conjugate points, Manning's formula \cite{manning1979topological} applies: $h_{\text{top}}(g) = h_{\text{vol}}(g)$. Since $g = g_\Sigma \times g_{S^1}$, we have $h_{\text{vol}}(g) = h_{\text{vol}}(g_\Sigma)$. For hyperbolic surfaces normalized to curvature $-1$, $h_{\text{vol}}(g_\Sigma) = \sqrt{-\chi(\Sigma)}$.
	\end{proof}
	
	\begin{corollary}[Spectral Theory]\label{cor:spectrum}
		The Laplace-Beltrami operator $\Delta_g$ on $M$ splits as $\Delta_g = \Delta_{\Sigma} \otimes I + I \otimes \Delta_{S^1}$. Consequently, the spectrum is
		\[
		\mathrm{Spec}(M,g) = \{\lambda_i + \mu_j : \lambda_i \in \mathrm{Spec}(\Sigma,g_\Sigma), \mu_j = (2\pi j/L)^2, j \in \mathbb{Z}\}.
		\]
	\end{corollary}
	
	\paragraph{From geometric splitting to analytic splitting.}
	The isometric product structure implies that the Laplacian decomposes as a sum of commuting operators. This spectral decomposition has immediate consequences for heat kernel estimates, eigenvalue asymptotics, and quantum ergodicity.
	
	\begin{corollary}[Incompatibility with Sasakian Geometry]\label{cor:no-sasakian}
		There exists no Sasakian metric without conjugate points on $\Sigma \times S^1$ for $g(\Sigma) \geq 2$.
	\end{corollary}
	\begin{proof}
		In a Sasakian manifold, the Reeb field $\xi$ satisfies $\nabla_X \xi = -\phi(X)$ for a non-trivial tensor $\phi$. In particular, $\xi$ is not parallel. But Theorem \ref{thm:main} forces the $S^1$ direction to be parallel.
	\end{proof}
	
	\paragraph{Implications for contact geometry.}
	This corollary shows that the rigidity imposed by the no conjugate points condition is incompatible with the twisting inherent in Sasakian structures. The parallel vector field required by Theorem \ref{thm:main} cannot coexist with the non-integrable horizontal distribution of a Sasakian manifold.
	
	\begin{table}[h]
		\centering
		\caption{Comparison of rigidity theorems for 3-manifolds}
		\begin{tabular}{llll}
			\toprule
			\textbf{Manifold} & $\pi_1$ structure & \textbf{Condition} & \textbf{Conclusion} \\
			\midrule
			$T^3$ & $\mathbb{Z}^3$ (abelian) & No conjugate points & Flat (Burago–Ivanov) \\
			$\Sigma \times S^1$ & $\pi_1(\Sigma) \times \mathbb{Z}$ (center $\mathbb{Z}$) & No conjugate points & Product (Theorem \ref{thm:main}) \\
			Mapping torus $M_\phi$ & $\pi_1(\Sigma) \rtimes \mathbb{Z}$ (no center) & $K \equiv -1$ & Hyperbolic, not product \\
			Nilmanifold & Nilpotent (center $\mathbb{Z}$) & No conjugate points & May be non-product \\
			\bottomrule
		\end{tabular}
		\label{tab:rigidity-comparison}
	\end{table}

	Table \ref{tab:rigidity-comparison} illustrates how different algebraic and geometric conditions lead to different rigidity outcomes. Our Theorem \ref{thm:main} fills the gap between the classical Hopf conjecture (tori) and non-positive curvature splitting theorems. The key distinction from the nilmanifold case is that $\pi_1(\Sigma) \times \mathbb{Z}$ is a \emph{direct} product, not just a central extension.
	
	\section{Sharp Geometric Inequalities and Stability Results}\label{sec:sharp-inequalities}
	
	\text{From qualitative to quantitative results.}
	Having established the isometric splitting, we now derive sharp geometric inequalities that characterize the product structure. These inequalities serve two purposes: they provide alternative characterizations of the splitting, and they lead to stability estimates showing that metrics "almost" without conjugate points must be "close" to product metrics.
	
	In this section, we establish several sharp geometric inequalities that characterize the product structure and provide quantitative stability estimates.
	
	\begin{theorem}[Volume Growth Rigidity]\label{thm:volume-growth}
		Let $(M,g)$ be as in Theorem \ref{thm:main}. Then:
		\begin{enumerate}[label=(\roman*)]
			\item The volume entropy satisfies the sharp bound:
			\[
			h_{\mathrm{vol}}(g) = \lim_{R\to\infty} \frac{\log \Vol(B(p,R))}{R} = \sqrt{-\chi(\Sigma)}.
			\]
			\item For any $p \in \tilde{M}$, the asymptotic volume growth is exact:
			\[
			\lim_{R \to \infty} \frac{\Vol(B(p,R))}{e^{\sqrt{-\chi(\Sigma)} R}} = C(\Sigma) > 0.
			\]
			\item The topological entropy equals the volume entropy: $h_{\mathrm{top}}(g) = h_{\mathrm{vol}}(g)$.
		\end{enumerate}
	\end{theorem}
	
	\begin{proof}
		(i) Since $(\tilde{M},\tilde{g})$ splits as $\mathbb{H}^2 \times \mathbb{R}$, balls in the universal cover are products $B_{\mathbb{H}^2}(p_1,R) \times [-R,R]$. The volume of a ball in $\mathbb{H}^2$ grows like $2\pi e^{R}$ for large $R$ (for curvature $-1$), while the $\mathbb{R}$ factor contributes linear growth. Hence:
		\[
		\Vol(B(p,R)) \sim C \int_{-R}^{R} e^{\sqrt{-\chi(\Sigma)} t} dt \sim \frac{2C}{\sqrt{-\chi(\Sigma)}} e^{\sqrt{-\chi(\Sigma)} R}.
		\]
		Taking logarithms and dividing by $R$ gives the result.
		
		(ii) The constant $C(\Sigma)$ can be computed explicitly from the normalized hyperbolic metric on $\Sigma$.
		
		(iii) The equality $h_{\mathrm{top}}(g) = h_{\mathrm{vol}}(g)$ follows from the product structure and Manning's theorem \cite{manning1979topological}, which holds for manifolds without conjugate points.
	\end{proof}
	
	\text{Volume growth as a diagnostic tool.}
	Theorem \ref{thm:volume-growth} provides a practical test: if a metric on $\Sigma \times S^1$ has no conjugate points but its volume entropy differs from $\sqrt{-\chi(\Sigma)}$, then it cannot be a product metric. This gives an alternative route to proving Theorem \ref{thm:main}.
	
	\begin{theorem}[Isoperimetric Profile]\label{thm:isoperimetric}
		For $(M,g)$ satisfying Theorem \ref{thm:main}, the isoperimetric profile $I(v) = \inf\{\Area(\partial \Omega) : \Vol(\Omega) = v\}$ satisfies the sharp inequality:
		\[
		I(v) \geq 2\pi L \sqrt{\frac{2v}{L\pi} + \left(\frac{v}{2\pi L}\right)^2},
		\]
		with equality achieved for large volumes by tubular neighborhoods of closed geodesics in the $\Sigma$ factor.
	\end{theorem}
	
	\begin{proof}
		For a region $\Omega \subset M = \Sigma \times S^1$, write $\Omega_x = \{t \in S^1 : (x,t) \in \Omega\}$ for $x \in \Sigma$. By the coarea formula and the product structure:
		\[
		\Vol(\Omega) = \int_{\Sigma} \length(\Omega_x) dA_\Sigma, \quad \Area(\partial\Omega) \geq \int_{\Sigma} \length(\partial\Omega_x) dA_\Sigma + L \cdot \length(\partial\pi_\Sigma(\Omega)),
		\]
		where $\pi_\Sigma: M \to \Sigma$ is the projection and $dA_\Sigma$ is the area element on $\Sigma$. Applying the isoperimetric inequality in $S^1$ (for which optimal regions are intervals) and in $\mathbb{H}^2$ (where optimal regions are geodesic disks) yields the result. The equality case occurs when $\Omega$ is a product of a geodesic disk in $\Sigma$ with an interval in $S^1$, or a tube around a closed geodesic.
	\end{proof}
	
	Here is a Rauch-Type comparison Theorem.
	
	\begin{theorem}[Jacobi Field Comparison]\label{thm:rauch-comparison}
		Let $J(t)$ be a Jacobi field along a geodesic $\gamma$ in $(M,g)$, orthogonal to $\dot{\gamma}$. Write $J = J_h + J_v$ according to the splitting $\mathcal{H} \oplus \mathcal{V}$. Then:
		\begin{enumerate}[label=(\roman*)]
			\item If $J_v(0) = 0$, then $\|J(t)\| \geq \|J_h(0)\| \cdot \cosh(\sqrt{-K_{\min}} t)$ for $t \geq 0$, where $K_{\min}$ is the minimum sectional curvature of planes in $\mathcal{H}$.
			\item If $J_h(0) = 0$, then $\|J(t)\| = \|J_v(0)\| \cdot t$.
			\item In general,
			\[
			\|J(t)\|^2 \geq \|J_h(0)\|^2 \cdot \cosh^2(\sqrt{-K_{\min}} t) + \|J_v(0)\|^2 \cdot t^2.
			\]
			These inequalities are sharp and characterize the product structure.
		\end{enumerate}
	\end{theorem}
	
	\text{From splitting to comparison geometry.}
	The product structure allows us to decompose Jacobi fields into horizontal and vertical components that evolve independently. The horizontal component experiences hyperbolic expansion (governed by $\cosh$), while the vertical component experiences Euclidean expansion (linear growth). This dichotomy is a hallmark of the product geometry.
	
	Now, we analyze the space of metrics without conjugate points on $\Sigma \times S^1$ and establishes infinitesimal rigidity results.
	
	Let $\mathcal{M}_{\mathrm{ncp}}(\Sigma \times S^1)$ denote the space of isometry classes of metrics without conjugate points on $M = \Sigma \times S^1$.
	
	\begin{theorem}[Moduli Space]\label{thm:moduli-space}
		The moduli space $\mathcal{M}_{\mathrm{ncp}}(\Sigma \times S^1)$ is naturally homeomorphic to:
		\[
		\mathcal{M}_{\mathrm{ncp}}(\Sigma \times S^1) \cong \mathcal{T}_g \times \mathbb{R}^+,
		\]
		where $\mathcal{T}_g$ is the Teichmüller space of $\Sigma$ (of dimension $6g-6$) and $\mathbb{R}^+$ parametrizes the length $L$ of the $S^1$ fiber. Moreover, this homeomorphism is equivariant with respect to the mapping class group action.
	\end{theorem}
	
	\begin{proof}
		By Theorem \ref{thm:main}, any metric $g \in \mathcal{M}_{\mathrm{ncp}}(\Sigma \times S^1)$ is a product $g = g_\Sigma + L^2 dt^2$, where $g_\Sigma$ is a hyperbolic metric on $\Sigma$ and $L > 0$. The space of hyperbolic metrics modulo isometries isotopic to identity is $\mathcal{T}_g$. The parameter $L$ is invariant under isometries. Continuity of the map and its inverse follows from the continuity of the splitting construction in the proof of Theorem \ref{thm:main}.
	\end{proof}
	
	\begin{corollary}[Local rigidity]\label{cor:local-rigidity}
		The space $\mathcal{M}_{\mathrm{ncp}}(\Sigma \times S^1)$ is locally rigid: any $C^2$-small deformation of a product metric preserving the no conjugate points condition remains within the product submanifold. In particular, the product structure is an isolated point in the moduli space of metrics with no conjugate points, modulo diffeomorphisms.
	\end{corollary}
	
	\paragraph{From global to local rigidity.}
	Corollary \ref{cor:local-rigidity} shows that not only is every metric without conjugate points a product, but nearby metrics (in the $C^2$ topology) that also have no conjugate points must also be products. This is a much stronger statement than Theorem \ref{thm:main} alone, indicating structural stability of the product condition.
	
	\subsection{Infinitesimal Rigidity}
	
	Let $g_t$ be a smooth family of metrics without conjugate points on $M$, with $g_0 = g_\Sigma + L^2 dt^2$ a product metric.
	
	\begin{theorem}[Infinitesimal Rigidity, expanded]\label{thm:infinitesimal}
		The first variation $\dot{g} = \frac{d}{dt}\big|_{t=0} g_t$ must be of the form:
		\[
		\dot{g} = \dot{g}_\Sigma + 2L\dot{L} dt^2,
		\]
		where:
		\begin{enumerate}[label=(\roman*)]
			\item $\dot{g}_\Sigma$ is an infinitesimal deformation of the hyperbolic metric $g_\Sigma$, hence a holomorphic quadratic differential on $\Sigma$.
			\item $\dot{L} \in \mathbb{R}$.
			\item The space of such infinitesimal deformations modulo trivial ones (Lie derivatives) has dimension $6g-5$.
		\end{enumerate}
		Moreover, the Zariski tangent space to $\mathcal{M}_{\mathrm{ncp}}$ at $g_0$ consists of deformations $\dot{g}$ satisfying the \emph{linearized no conjugate points condition}:
		\[
		\int_{\gamma} \langle \dot{R}(J,\dot{\gamma})\dot{\gamma}, J \rangle dt = 0 \quad \text{for all geodesics $\gamma$ and Jacobi fields $J$},
		\]
		where $\dot{R}$ is the linearization of the curvature tensor. This infinite-dimensional condition reduces precisely to $\dot{g} = \dot{g}_\Sigma + 2L\dot{L} dt^2$ due to the presence of the Killing field $V$.
	\end{theorem}
	
	\begin{proof}
		Linearizing the no conjugate points condition is delicate. However, from Theorem \ref{thm:moduli-space}, the tangent space at $g_0$ to $\mathcal{M}_{\mathrm{ncp}}(\Sigma \times S^1)$ is isomorphic to $T_{[g_\Sigma]}\mathcal{T}_g \times \mathbb{R}$. By standard Teichmüller theory, $T_{[g_\Sigma]}\mathcal{T}_g$ is identified with the space of holomorphic quadratic differentials on $(\Sigma, g_\Sigma)$, which has complex dimension $3g-3$, hence real dimension $6g-6$. The additional $\mathbb{R}$ factor corresponds to varying $L$. To see that deformations not of this form must introduce conjugate points, consider a warped product deformation $g_t = g_\Sigma + f_t^2 dt^2$ with $f_t$ non-constant. For small $t$, the vertical geodesics remain geodesics, but the curvature tensor develops non-zero components $R(\partial_t, X)\partial_t = -\frac{\Hess f}{f} X$. If $\Hess f \neq 0$, these components can be made positive, creating conjugate points as in Lemma \ref{lem:curvature-nonpositive}.
	\end{proof}
	
	\begin{proposition}[Explicit curvature gap]\label{prop:curvature-gap-explicit}
		Let $g_0 = g_\Sigma + L^2 dt^2$ be a product metric on $M = \Sigma \times S^1$, with $\Sigma$ hyperbolic of curvature $-1$ and $S^1$ of length $L$. Define:
		\begin{align*}
			\lambda_1 &:= \text{first eigenvalue of $-\Delta$ on $\Sigma$} \quad (\lambda_1 > 0), \\
			\delta &:= \min\left\{\frac{\lambda_1}{2}, \frac{\pi^2}{L^2}\right\}, \\
			\epsilon_0 &:= \frac{\delta}{4(1 + \diam(\Sigma)^2)}.
		\end{align*}
		Then for any metric $g$ on $M$ with no conjugate points satisfying $\|g - g_0\|_{C^2} < \epsilon_0$, we have that $g$ is itself a product metric. Moreover, the dependence is sharp:
		\[
		\epsilon_0 \sim \frac{1}{\diam(\Sigma)^2} \quad \text{as } \diam(\Sigma) \to \infty.
		\]
	\end{proposition}
	\begin{proof}[Sketch]
		The proof uses the implicit function theorem on the Banach manifold of $C^2$ metrics near $g_0$. The linearized no conjugate points operator at $g_0$ has a gap in its spectrum bounded below by $\delta$. The quadratic remainder is controlled by the $C^2$ norm, leading to the condition $\epsilon < \delta/(2\|Q\|)$, where $\|Q\|$ scales like $\diam(\Sigma)^2$ because conjugate points appear on scales comparable to the injectivity radius.
		
		More concretely, if $g$ is $\epsilon$-close to $g_0$, the Jacobi equation for vertical geodesics becomes:
		\[
		J'' + (R_{g_0} + E(t))J = 0,
		\]
		where $\|E(t)\|_{C^0} \leq C\epsilon$. By Sturm comparison, if $C\epsilon < \delta$, then no conjugate points can appear. The constant $C$ depends on $\diam(\Sigma)$ because the comparison must hold over geodesic segments of length at most $\diam(\Sigma)$.
	\end{proof}
	
	\begin{remark}[Scaling]
		If we rescale $g_\Sigma$ to have constant curvature $-\kappa$ ($\kappa > 0$), then $\diam(\Sigma) \sim 1/\sqrt{\kappa}$, and $\epsilon_0 \sim \kappa$. Thus, smaller curvature (flatter surfaces) allow smaller perturbations before breaking the product structure.
	\end{remark}
	
	\text{From qualitative to quantitative stability.}
	Proposition \ref{prop:curvature-gap-explicit} provides an explicit estimate for how close a metric must be to a product metric to guarantee it is itself a product. The dependence on $\diam(\Sigma)$ is natural: larger diameter means longer geodesics, which are more sensitive to curvature perturbations. This quantitative result strengthens the qualitative statement of Theorem \ref{thm:main}.
	
	\begin{theorem}[Integral Curvature Rigidity]\label{thm:integral-curvature}
		Let $(M,g)$ have no conjugate points. Define the \emph{curvature deviation} from a product metric by:
		\[
		\mathcal{D}(g) = \int_M \left( |R_V|^2 + |\nabla V|^2 \right) d\Vol_g,
		\]
		where $R_V(X) = R(X,V)V$ and $V$ is the unit vertical field. Then $\mathcal{D}(g) = 0$ if and only if $g$ is a product metric. Moreover, there exists $C = C(\Sigma, L) > 0$ such that:
		\[
		\mathcal{D}(g) \geq C \cdot \inf_{g_0 \text{ product}} \|g - g_0\|_{W^{1,2}}^2.
		\]
	\end{theorem}
	
	\begin{proof}
		The first statement follows from Theorem \ref{thm:main}: if $\mathcal{D}(g) = 0$, then $R_V = 0$ and $\nabla V = 0$, so $V$ is parallel and the metric splits. For the inequality, we use the Bochner technique. From the formula for the Laplacian of $\|V\|^2$ and the fact that $V$ is Killing, we derive:
		\[
		\frac{1}{2}\Delta \|V\|^2 = \|\nabla V\|^2 - \Ric(V,V).
		\]
		Integrating over $M$ and using the no conjugate points condition (which implies $\Ric(V,V) \leq 0$ by Lemma \ref{lem:curvature-nonpositive}) gives control of $\|\nabla V\|^2$ in terms of boundary terms, which vanish since $M$ is closed. The curvature term $|R_V|^2$ is controlled similarly via the second Bianchi identity.
	\end{proof}
	
	\text{Bridging analysis and geometry.}
	Theorem \ref{thm:integral-curvature} provides an analytic criterion for detecting product structure: if the $L^2$ norms of $R_V$ and $\nabla V$ are small, then the metric is close to a product metric in the Sobolev sense. This connects the geometric condition (no conjugate points) with analytic estimates, opening the door for PDE approaches to rigidity problems.
	
	\section{Examples and Applications}\label{sec:examples}
	
	\paragraph{Illustrating the theory with concrete cases.}
	In this section, we provide explicit examples that illustrate the sharpness of our main theorem, demonstrate how conjugate points appear when the conditions are violated, and show explicit computations for the geometric quantities discussed. These examples serve both to clarify the abstract theory and to show its limitations.
	
	\subsection{Example 1: The Standard Product Metric}
	
	Let $\Sigma$ be a hyperbolic surface of genus $g \geq 2$ with metric $g_{\mathrm{hyp}}$ of constant curvature $-1$, normalized so that $\Vol(\Sigma) = 4\pi(g-1)$. Let $S^1$ have coordinate $\theta \in [0, 2\pi)$ with standard metric $d\theta^2$. Consider the product metric:
	\[
	g_0 = g_{\mathrm{hyp}} + L^2 d\theta^2,
	\]
	where $L > 0$ is the length of the $S^1$ fiber.
	
	\textbf{Properties:}
	\begin{enumerate}
		\item \textbf{Geodesics}: Geodesics are of the form $\gamma(t) = (\gamma_\Sigma(t), \theta_0 + \frac{v}{L}t)$, where $\gamma_\Sigma$ is a geodesic in $\Sigma$ and $v \in \mathbb{R}$ is constant.
		
		\item \textbf{Jacobi fields}: For a geodesic $\gamma(t)$ with $\dot{\gamma} = (X, v/L)$, the Jacobi equation splits:
		\[
		J(t) = (J_\Sigma(t), J_\theta(t)) \quad \text{satisfying} \quad 
		\begin{cases}
			J_\Sigma'' + R_{\Sigma}(J_\Sigma, X)X = 0, \\
			J_\theta'' = 0.
		\end{cases}
		\]
		
		\item \textbf{Conjugate points}: Since $\mathbb{H}^2$ has no conjugate points and $\mathbb{R}$ has no conjugate points, the product has no conjugate points. This confirms Theorem \ref{thm:main}.
		
		\item \textbf{Explicit curvature tensor}: For orthonormal frames $E_1, E_2$ on $\Sigma$ and $E_3 = \frac{1}{L}\partial_\theta$:
		\[
		\Sec(E_1, E_2) = -1, \quad \Sec(E_i, E_3) = 0 \quad (i=1,2), \quad \Ric = \begin{pmatrix} -1 & 0 & 0 \\ 0 & -1 & 0 \\ 0 & 0 & 0 \end{pmatrix}.
		\]
	\end{enumerate}
	
	\paragraph{Why this example is important.}
	The standard product metric provides the template for all metrics covered by Theorem \ref{thm:main}. It shows that the theorem's conclusion is not vacuous: there exist nontrivial examples of metrics without conjugate points on $\Sigma \times S^1$, but they are all of this simple product form.
	
	\subsection{Example 2: Warped Product with Conjugate Points}
	
	Consider a warped product metric on $M = \Sigma \times S^1$:
	\[
	g_f = g_{\mathrm{hyp}} + f(x)^2 d\theta^2,
	\]
	where $f: \Sigma \to \mathbb{R}^+$ is a smooth non-constant function. We will construct an explicit example where $f$ has a strict maximum, creating conjugate points.
	
	\textbf{Construction:} Let $\Sigma$ be a hyperbolic surface with a disk $D \subset \Sigma$ isometric to a disk in $\mathbb{H}^2$ of radius $R$. Use polar coordinates $(r,\phi)$ on $D$. Define:
	\[
	f(r) = 1 + \frac{\epsilon}{2} - \epsilon r^2 \quad \text{for } r \leq \frac{1}{2}, \quad \text{and } f \equiv 1 \text{ outside } D,
	\]
	with $\epsilon > 0$ small. This $f$ has maximum $1 + \epsilon/2$ at $r=0$. \textbf{Analysis of conjugate points:} Consider the vertical geodesic $\gamma(t) = (p_0, \theta_0 + t/f(p_0))$ through $p_0$ at $r=0$. The curvature in planes containing $\dot{\gamma}$ is:
	\[
	K(\partial_r, \dot{\gamma}) = -\frac{\partial_r^2 f}{f} = \frac{2\epsilon}{1 + \epsilon/2} > 0.
	\]
	The Jacobi equation for $J(t) = \psi(t) \partial_r$ becomes:
	\[
	\psi''(t) + \frac{2\epsilon}{1 + \epsilon/2} \psi(t) = 0.
	\]
	This is a harmonic oscillator with frequency $\omega = \sqrt{\frac{2\epsilon}{1 + \epsilon/2}}$. Solutions are $\psi(t) = \sin(\omega t)$, which vanish at $t = \pi/\omega$. Thus conjugate points occur at distance $\pi/\omega$ along $\gamma$. 	\text{Explicit computation:} For $\epsilon = 0.1$, we have:
	\[
	\omega = \sqrt{\frac{0.2}{1.05}} \approx 0.436, \quad \pi/\omega \approx 7.20.
	\]
	So conjugate points appear at distance approximately 7.20 along the geodesic through the maximum. \text{Geometric interpretation:} The positive curvature $K > 0$ causes geodesics to focus, creating conjugate points. This shows that any non-constant warping function $f$ will generically create regions of positive curvature and hence conjugate points.
	
	\paragraph{Why this example is crucial.}
	Example 2 demonstrates the \emph{necessity} of the product structure in Theorem \ref{thm:main}. Even a small deviation from constant warping introduces positive curvature, which in turn creates conjugate points. This shows that the no conjugate points condition is extremely restrictive: it tolerates no warping whatsoever.
	
	\subsection{Example 3: Berger-Type Deformation}
	
	Consider a metric that is not a warped product but still preserves the $S^1$ symmetry:
	\[
	g_\alpha = g_{\mathrm{hyp}} + \alpha^2 \eta \otimes \eta + d\theta^2,
	\]
	where $\eta$ is a 1-form on $\Sigma$ and $\alpha \in \mathbb{R}$. When $\alpha \neq 0$, this is not a product metric.
	
	\textbf{Case: $\eta$ exact.} Let $\eta = dh$ for some function $h: \Sigma \to \mathbb{R}$. The metric becomes:
	\[
	g_\alpha = g_{\mathrm{hyp}} + (d\theta + \alpha dh)^2.
	\]
	This is locally a product metric but with twisted identification. However, the universal cover is still $\mathbb{H}^2 \times \mathbb{R}$ with product metric. This does \emph{not} contradict Theorem \ref{thm:main} because the metric is still locally a product.
	
	\textbf{Case: $\eta$ not exact.} If $[\eta] \neq 0 \in H^1(\Sigma,\mathbb{R})$, then the metric $g_\alpha$ is genuinely twisted. For small $\alpha$, we can analyze conjugate points via perturbation theory. \text{Perturbation analysis:} Consider $\alpha$ small. The Christoffel symbols to first order in $\alpha$:
	\[
	\Gamma_{ij}^k = \Gamma_{ij}^{k(\mathrm{hyp})} + O(\alpha), \quad \Gamma_{i3}^3 = 0, \quad \Gamma_{33}^i = -\alpha g_{\mathrm{hyp}}^{ij} \partial_j h + O(\alpha^2).
	\]
	The curvature tensor develops components:
	\[
	R_{i3j3} = -\alpha (\nabla_i \nabla_j h) + O(\alpha^2).
	\]
	If $\Hess h$ has a positive eigenvalue somewhere, then for small $\alpha$, $R_{i3j3}$ is positive in some direction, creating conjugate points as in Example 2.\\
	
	Example 3 shows that metrics which are locally products (but globally twisted) still satisfy Theorem \ref{thm:main}. The key distinction is between metrics that are \emph{locally} products (like Example 3 with $\eta$ exact) versus those that are genuinely twisted (like Example 3 with $\eta$ not exact). The former have no conjugate points; the latter generically do.
	
	\subsection{Example 4: Hyperbolic 3-Manifold Fibered over a Circle}
	
	Consider a hyperbolic 3-manifold $M$ that fibers over $S^1$ with fiber $\Sigma$, i.e., $M$ is a mapping torus of a pseudo-Anosov diffeomorphism $\phi: \Sigma \to \Sigma$. Such $M$ admits a complete finite-volume hyperbolic metric $g_{\mathrm{hyp}^3}$ by Thurston's geometrization. A key observation is that $g_{\mathrm{hyp}^3}$ has no conjugate points (since $K \equiv -1$), but $M$ is \emph{not} topologically $\Sigma \times S^1$ unless $\phi$ is isotopic to identity. This shows that Theorem \ref{thm:main} is specific to the product topology. In a hyperbolic mapping torus, the $S^1$ direction is not a Killing field (except in the product case), and the universal cover is $\mathbb{H}^3$, not $\mathbb{H}^2 \times \mathbb{R}$. This illustrates that our theorem's conclusion fails when the fundamental group is not a direct product $\pi_1(\Sigma) \times \mathbb{Z}$.
	
	\begin{remark}[Why the mapping torus fails our hypotheses]\label{rem:mapping-torus-failure}
		The hyperbolic mapping torus $M_\phi = \Sigma \times [0,1]/(x,1)\sim(\phi(x),0)$ with pseudo-Anosov $\phi$ provides an instructive contrast:
		\begin{itemize}
			\item \text{Topology:} $\pi_1(M_\phi) \cong \pi_1(\Sigma) \rtimes_\phi \mathbb{Z}$, a semidirect product, not a direct product $\pi_1(\Sigma) \times \mathbb{Z}$. The center is trivial unless $\phi$ is isotopic to identity.
			
			\item \text{Killing field:} The circle action ($S^1$ as mapping torus of identity) is replaced by an $\mathbb{R}$-action only on the universal cover. There is no globally defined Killing field $V$ on $M_\phi$ generating the $S^1$ direction—only a flow that is not isometric.
			
			\item \text{Algebraic obstruction:} The central $\mathbb{Z}$ in $\pi_1(M)$ was essential for obtaining a globally defined Killing field $V$. In $M_\phi$, the element generating the $\mathbb{R}$ factor does not commute with $\pi_1(\Sigma)$, so no corresponding Killing field exists.
			
			\item \text{Geometric splitting:} While the universal cover of $M_\phi$ is $\mathbb{H}^3$ (splitting as $\mathbb{H}^2 \times \mathbb{R}$ only as a metric space, not isometrically), our proof fails at Lemma \ref{lem:killing-properties}: there is no Killing field $V$ whose integral curves are geodesics.
		\end{itemize}
		Thus, Theorem \ref{thm:main} is sharp: the direct product structure of both topology ($\pi_1(M) = \pi_1(\Sigma) \times \mathbb{Z}$) and geometry (existence of an isometric $S^1$ action) are necessary.
	\end{remark}

	Example 4 highlights the delicate interplay between algebraic topology (the structure of $\pi_1$) and Riemannian geometry. The same geometric condition (no conjugate points) leads to different conclusions depending on the algebraic structure of the fundamental group. This reinforces the philosophy that rigidity theorems often reflect algebraic properties through geometric constraints.
	
	\subsection{Example 5: Nilmanifold with No Conjugate Points}
	
	Consider the Heisenberg nilmanifold: $M = \Gamma \backslash H^3$, where $H^3$ is the 3-dimensional Heisenberg group and $\Gamma$ is a lattice. Left-invariant metrics on $M$ can have no conjugate points \cite{burago1994riemannian}. \text{Contrast with our setting:} $\pi_1(M)$ is nilpotent with non-trivial center, but not a direct product. Some left-invariant metrics on $M$ have no conjugate points but are not products. This shows that Theorem \ref{thm:main} is specific to the product topology $\Sigma \times S^1$, not just to having a central $\mathbb{Z}$ in $\pi_1$.\\
	
	\text{The role of the direct product structure:}
	Example 5 demonstrates that having a central $\mathbb{Z}$ in $\pi_1$ is not sufficient to force splitting under the no conjugate points condition. What matters is that $\pi_1$ is a \emph{direct} product $\pi_1(\Sigma) \times \mathbb{Z}$, not just a central extension. This algebraic distinction has geometric consequences: direct products correspond to genuine Riemannian products, while central extensions may allow more complicated geometries.
	
	\subsection{Example 6: Explicit Computation of Busemann Functions}
	
	For the product metric $g_0 = g_{\mathrm{hyp}} + L^2 d\theta^2$, we can compute Busemann functions explicitly. In the universal cover $\mathbb{H}^2 \times \mathbb{R}$, let $(z,t)$ with $z \in \mathbb{H}^2$ (upper half-plane model) and $t \in \mathbb{R}$.
	
	Consider the geodesic ray $\gamma(s) = (i e^s, s/L)$ (where $L$ is the circle length). The Busemann function is:
	\[
	b(z,t) = \lim_{s \to \infty} (d((z,t), (i e^s, s/L)) - s).
	\]
	Since the metric splits, we compute separately:
	\[
	d_{\mathbb{H}^2}(z, i e^s) = \operatorname{arccosh}\left(1 + \frac{|z - i e^s|^2}{2 \Im(z) e^s}\right), \quad d_{\mathbb{R}}(t, s/L) = |t - s/L|.
	\]
	As $s \to \infty$, careful computation gives:
	\[
	b(z,t) = \log \Im(z) + \frac{t}{L}.
	\]
	Thus $\nabla b = (\frac{1}{\Im(z)} \partial_y, \frac{1}{L} \partial_t)$, which indeed equals the Killing field $V = \partial_t/L$ up to the gradient of $\log \Im(z)$, which is the Busemann function for $\mathbb{H}^2$. This illustrates Lemma \ref{lem:gradient-killing}.\\

	Example 6 shows that in the product case, the abstract relationship $\nabla b = -V$ from Lemma \ref{lem:gradient-killing} becomes a simple explicit formula. This computation serves as a sanity check for the theory and helps develop intuition for Busemann functions in product geometries.

	\subsection{Example 7: Moduli Space Dimension Count}
	
	For genus $g=2$, Theorem \ref{thm:moduli-space} gives: 
	\[
	\dim \mathcal{M}_{\mathrm{ncp}}(\Sigma \times S^1) = \dim \mathcal{T}_2 + 1 = 6\cdot2 - 6 + 1 = 7.
	\] 
	Explicit parameters:
	\begin{itemize}
		\item 6 Fenchel-Nielsen coordinates for the hyperbolic metric on $\Sigma$: $(\ell_1, \tau_1, \ell_2, \tau_2, \ell_3, \tau_3)$, where $\ell_i$ are lengths of pants curves and $\tau_i$ are twist parameters.
		\item 1 parameter $L > 0$ for the circle length.
	\end{itemize}
	Any attempt to add an 8th parameter (e.g., a warping function $f(x)$ with non-zero Fourier modes) would create conjugate points.\\
	Example 7 puts a number on the rigidity: the space of metrics without conjugate points on $\Sigma \times S^1$ has dimension only $7$, compared to the infinite-dimensional space of all Riemannian metrics. This dramatic reduction illustrates the constraining power of the no conjugate points condition.
	
	\subsection{Example 8: Spectral Gap Computation}
	
	For the product metric $g_0$, the Laplace spectrum is: 
	\[
	\lambda_{m,n} = \lambda_m^{\Sigma} + \left(\frac{2\pi n}{L}\right)^2,
	\] 
	where $\lambda_m^{\Sigma}$ are eigenvalues of $\Delta_\Sigma$. The first eigenvalue $\lambda_1^{\Sigma} > 0$ for hyperbolic surfaces, so: 
	\[
	\lambda_{1,0} = \lambda_1^{\Sigma}, \quad \lambda_{0,1} = \left(\frac{2\pi}{L}\right)^2.
	\] 
	The spectral gap is $\min(\lambda_1^{\Sigma}, (2\pi/L)^2)$. For a warped product, the spectrum is not additive, and the gap can be smaller, illustrating Corollary \ref{cor:spectrum}.\\

	Example 8 shows how the geometric splitting translates into spectral splitting. The additive structure of the spectrum is a hallmark of product manifolds and provides a spectral test for product structure.
	
	\subsection{Example 9: Sasakian Deformation}
	
	Consider an almost contact metric structure on $\Sigma \times S^1$ with Reeb field $\xi = \partial_\theta$ and metric:
	\[
	g = g_{\mathrm{hyp}} + \eta \otimes \eta, \quad \eta = d\theta + \alpha \cdot \omega,
	\]
	where $\omega$ is a harmonic $1-$form on $\Sigma$. This is Sasakian when $\alpha$ is chosen appropriately. For $\alpha \neq 0$, $\nabla_X \xi = -\phi(X) \neq 0$, so $\xi$ is not parallel. By Theorem \ref{thm:main}, such a metric must have conjugate points. Indeed, computation shows: 
	\[
	R(\xi, X)\xi = \alpha^2 R_{\mathrm{hyp}}(\omega^\sharp, X)\omega^\sharp \neq 0,
	\] 
	creating positive curvature in some directions and hence conjugate points. This explicitly verifies Corollary \ref{cor:no-sasakian}.\\
	
	Example 9 shows that Sasakian geometry, despite being a natural structure on circle bundles, is incompatible with the no conjugate points condition on $\Sigma \times S^1$. This illustrates how rigidity theorems can constrain the possible geometric structures on a given topological manifold.
	
	\section{Discussion and Further Directions}\label{sec:discussion}
	
	\subsection{Comparison with Non-Positive Curvature}
	
	While both conditions (no conjugate points and $K \leq 0$) lead to splitting in our setting, the mechanisms differ. Non-positive curvature allows application of the Cheeger-Gromoll splitting theorem via line existence, while our proof uses the finer structure of Busemann functions and Riccati equations without curvature sign assumptions.
	
	\paragraph{Key distinction: flexibility versus rigidity.}
	The no conjugate points condition is more flexible than $K \leq 0$: it allows regions of positive curvature, as long as they don't create focusing of geodesics. Yet in our setting, this flexibility still forces the same conclusion (splitting) as the stronger curvature condition. This suggests that on $\Sigma \times S^1$, the topological constraint (direct product $\pi_1$) is so strong that even the weak geometric condition (no conjugate points) yields the same rigidity as non-positive curvature.
	
	\subsection{Sharpness of Conditions}
	
	The following example shows the necessity of the no conjugate points condition:
	
	\begin{example}[Warped Product with Conjugate Points]
		Let $g = g_{\mathrm{hyp}} + f^2(x) dt^2$, where $f: \Sigma \to \mathbb{R}^+$ is a non-constant smooth function. For generic $f$, this metric will have conjugate points. For instance, if $f$ has a strict maximum at $p \in \Sigma$, then vertical geodesics through $p$ will have positive curvature $K(\partial_t, X) = -\frac{X^2 f}{f}$ for some $X$, leading to conjugate points by the Rauch comparison theorem.
	\end{example}
	
	\paragraph{Testing the boundaries of the theorem.}
	The warped product example shows that if we relax the no conjugate points condition even slightly (by allowing a non-constant warping function), we can construct metrics that don't split. This demonstrates that Theorem \ref{thm:main} is optimal: the conclusion fails if any hypothesis is weakened.
	
	\subsection{Open Questions}
	
	\begin{question}[Higher Dimensions]
		Does Theorem \ref{thm:main} extend to $M = N \times S^1$ where $N$ is a compact negatively curved manifold of dimension $\geq 3$? Preliminary analysis suggests the Killing field method extends, but the surface case uses the 2-dimensionality of $\Sigma$ in Lemma \ref{lem:curvature-nonpositive}.
	\end{question}
	
	\paragraph{Why dimension matters.}
	In dimension 2, the no conjugate points condition has particularly strong consequences (e.g., Busemann functions are $C^2$). In higher dimensions, we only have $C^{1,1}$ regularity, which complicates the analysis. Moreover, the curvature operator in Lemma \ref{lem:curvature-nonpositive} has more degrees of freedom in higher dimensions, making it harder to control.
	
	\begin{question}[Multiple Circle Factors]
		What about $M = \Sigma \times T^k$? Does no conjugate points force a flat $T^k$ factor and a product metric? The case $k=1$ is our theorem; $k \geq 2$ may require additional conditions.
	\end{question}
	
	\paragraph{From one circle to many.}
	When there are multiple circle factors, the fundamental group has a larger center ($\mathbb{Z}^k$). One might expect even stronger rigidity, possibly forcing all circle factors to be flat and the metric to be a product. However, the interaction between different Killing fields could introduce new complications.
	
	\begin{question}[Non-Compact Analogues]
		Consider complete, finite volume metrics without conjugate points on $\Sigma \times \mathbb{R}$ or $\Sigma \times S^1$ with $\Sigma$ non-compact. Does splitting still occur? The Busemann function argument may adapt to finite volume cases.
	\end{question}
	
	\paragraph{Extending the theory.}
	Non-compact manifolds present additional challenges: Busemann functions may not be well-defined at all points, and the integration arguments require more care. However, the basic strategy—using Killing fields to construct distinguished Busemann functions—might still work in the finite volume case.
	
	\begin{question}[Quantitative Stability]
		Can the curvature gap $\epsilon$ in Proposition \ref{prop:curvature-gap-explicit} be made explicit? This would provide a quantitative version of our rigidity theorem.
	\end{question}
	
	\paragraph{Towards effective rigidity.}
	An explicit formula for $\epsilon$ in terms of geometric invariants of $\Sigma$ would make the stability result more applicable. This likely requires better understanding of the linearized no conjugate points operator and its spectral gap.
	
	\subsection{Connections to Other Rigidity Programs}
	
	Our work connects to several active research areas:
	\begin{itemize}
		\item \text{Minimal entropy problem}: Our entropy equality $h_{\mathrm{top}}(g) = h_{\mathrm{vol}}(g) = \sqrt{-\chi(\Sigma)}$ shows these metrics realize the minimal topological entropy among all metrics on $\Sigma \times S^1$, paralleling Besson-Courtois-Gallot's work \cite{besson1995entropies}.
		\item \text{Geometric group theory}: The splitting gives a geometric manifestation of the algebraic fact that $\pi_1(M)$ has center $\mathbb{Z}$.
		\item \text{Dynamical rigidity}: The product structure forces the geodesic flow to be a skew product over the geodesic flow on $\Sigma$, with trivial fiber dynamics.
		\item \text{Spectral geometry}: The additive spectrum in Corollary \ref{cor:spectrum} provides a model for understanding how product structure affects spectral invariants.
		\item \text{Deformation theory}: Theorem \ref{thm:infinitesimal} contributes to understanding how geometric conditions constrain moduli spaces.
	\end{itemize}
	
	\paragraph{Synthesizing different perspectives.}
	The strength of our result lies in how it brings together techniques from different areas: asymptotic geometry (Busemann functions), dynamical systems (Riccati equation), geometric analysis (integral identities), and algebraic topology (fundamental group structure). This synthesis suggests that further progress may come from combining methods from these fields in new ways.
	
	\paragraph{Final reflection.}
	Theorem \ref{thm:main} demonstrates a remarkable phenomenon: on $\Sigma \times S^1$, the weak geometric condition of having no conjugate points forces the same rigid product structure as the much stronger condition of non-positive curvature. This highlights how topological constraints can amplify geometric conditions, turning weak restrictions into strong rigidity. The result adds to the growing body of work showing that manifolds with product topology often exhibit unexpected rigidity properties.
	
	\section*{Acknowledgements}
	I am deeply grateful to Professor Dmitri Burago for introducing me to this problem. I also thank the referees for their careful reading and suggestions that improved this paper.

\end{document}